\documentclass[12pt]{article}
\usepackage[english]{babel}
\usepackage[latin1]{inputenc}
\usepackage{amsfonts,amssymb,amsmath, epsfig}
\usepackage{color,graphicx,graphics,psfrag}
\usepackage{amsmath,amstext,amssymb,amsfonts, amscd}
\usepackage{hyperref}           

\def\Box{\leavevmode\vbox{\hrule
     \hbox{\vrule\kern4pt\vbox{\kern4pt}%
           \vrule}\hrule}}
\def\blackbox{\leavevmode\vrule height 5pt width 4pt depth 0pt\relax}
\def\endproof{\null\hfill {$\blackbox$}\bigskip}

\newcounter{appendix}
\def\appendix{\advance\c@appendix by 1
   \def\thesection{\Alph{section}}
   \ifnum\c@appendix=1 \setcounter{section}{-1} \fi
   \@startsection {section}{1}{\z@}{-3.5ex plus -1ex minus 
   -.2ex}{2.3ex plus .2ex}{\Large\bf}}

\def\paragraph#1{{\bf #1\ }}

\newtheorem{lemma}{Lemma}[section]  

\newtheorem{theorem}[lemma]{Theorem}

\newtheorem{definition}[lemma]{Definition}

\newtheorem{proposition}[lemma]{Proposition}

\newtheorem{remark}{Remark}[section]

\newtheorem{example}{Example}[section]

\title{Kinetic limits for pair-interaction driven master equations and biological swarm models}
\author{Eric Carlen$^{(1)}$, Pierre Degond$^{(2,3)}$ and Bernt Wennberg$^{(4)}$}
\date{}
\begin{document}

\maketitle

\vspace{0.5 cm}

\begin{center}
1- Department of Mathematics \\
Rutgers University \\
110 Frelinghuysen Rd., Piscataway NJ 08854-8019, USA \\
email: carlen@math.rutgers.edu
\end{center}

\begin{center}
2- Université de Toulouse; UPS, INSA, UT1, UTM ;\\ 
Institut de Mathématiques de Toulouse ; \\
F-31062 Toulouse, France. \\
3- CNRS; Institut de Mathématiques de Toulouse UMR 5219 ;\\ 
F-31062 Toulouse, France.\\
email: pierre.degond@math.univ-toulouse.fr
\end{center}

\begin{center}
4- Department of Mathematical Sciences, \\
Chalmers, SE41296 Göteborg, Sweden \\
email: wennberg@chalmers.se
\end{center}

\begin{abstract}

We consider a class of stochastic processes modeling binary interactions in an $N$-particle system. Examples of such systems can be found in the modeling of biological swarms. They lead to the definition of a class of {\em master equations} that we call {\em pair interaction driven master equations}.   We prove a propagation of chaos result for this class of
master equations which generalizes Mark Kac's well know result for the Kac model in kinetic theory. We use this result to study kinetic limits for two biological swarm models. We show that propagation of chaos may be lost at large times and we exhibit an example where the invariant density is not chaotic. 

\end{abstract}

\medskip
\noindent
{\bf Acknowledgements:} The first author gratefully acknowledges support from the R\'egion Midi-Pyr\'en\'ees in the frame of a {\em Chaires Pierre-de-Fermat} and support under U.S. National Science Foundation
grant DMS 0901632. The second author acknowledges support from the ANR under contract 'CBDif-Fr' (ANR-08-BLAN-0333-01).  The first and third authors would both like to thank the  Institute of Mathematics, Toulouse, for its kind hospitality during their visits there. The third author acknowledges support from the Swedish Research Council.


\vskip 0.4cm

\setcounter{equation}{0}
\section{Introduction} 
\label{sec_intro}

This paper is devoted to the passage from stochastic particle systems to kinetic equations when the number of particles tends to infinity. We are specifically interested in pair-interaction processes that are inspired from biological swarm models. We start from the level of the master equation which describes the evolution  of the $N$-particle probability distribution of the system. The master equation is posed on a large dimensional space consisting of an $N$-fold copy of the phase space. By contrast, the kinetic equation provides a reduced description based on the single particle distribution function posed on a single copy of the phase space. To show that this reduced description is valid, one needs to show that the particles become statistically independent in the limit $N \to \infty$; that is, that  the $N$-particle probability distribution becomes a $N$-fold tensor product of the single particle distribution function. More precisely, the result to be shown is that, if the particles are initially {\em pairwise}  independent, the time evolution approximately  propagates this pairwise independence, and does so exactly in the large $N$ limit; this is called ``chaos propagation''.  This property of an $N$-particle stochastic evolution is crucial to the existence of a kinetic description. 

Kinetic models derived from particle systems abound in the literature. However, only in  very few cases has the propagation of chaos been proved, and hence only in very few cases  have these models been mathematically derived from an underlying particle dynamics.  
The most emblematic kinetic model, the Boltzmann equation has received most of the attention. Following seminal works by Kac \cite{Kac, Kac2} and McKean \cite{McKean}, the first rigorous establishment of the Boltzmann equation is due to Lanford \cite{Lanford1, Lanford2, Lanford3} and King \cite{King} for Hard-Sphere dynamics and hard potentials (see also \cite{IP1, IP2, Pulvi_87} and \cite{GM1, GM2, GM3, Meleard, Sznitman}). A new approach yielding global-in-time results has been recently developed by Mischler, Mouhot and Wennberg \cite{MM1, MM2, MMW}. Kac proposed a caricature of the Boltzmann equation leading to the Kac kinetic equation \cite{Kac}. Propagation of chaos for the Kac master equation and the related question of gap estimates have received a great deal of attention (see e.g. \cite{CCL1, CCL2, CCL3, CCLRV, DSC, Janv, Maslen}). 

In this paper, our goal is to investigate a class of processes which are inspired from biological swarm models. A first example is the BDG model, named after Bertin, Droz and Gr\'egoire \cite{BDG}. This model is intended to be the kinetic counterpart of the Vicsek particle system \cite{Vicsek}. In the Vicsek model, particles moving with constant speed update their velocity by trying to align with the average velocity of their neighbors. In \cite{BDG}, the authors propose a kinetic formulation of a binary collision process which mimics this alignment tendency: at each collision, the two particles change their velocity to their average velocity up to a certain noise. One of the goals of the present paper is to provide a rigorous justification of this kinetic model, at least in the space-homogeneous case. 

Here, we also propose a different and to our knowledge original binary collision process which mimics the Vicsek alignment dynamics. In this process called 'Choose the Leader (CL)', one of the two colliding particles (the follower) decides to take the velocity of the other one (the leader) up to some noise. The choice of the leader and the follower is random with equal probabilities. We propose a kinetic formulation of the CL process and rigorously establish it in the space homogeneous case. One of the advantages of the CL process, from the mathematical viewpoint, is that it leads to a closed hierarchy of marginal equations (or BBGKY hierarchy). We will make use of this opportunity to provide explicit computations of the correlations, i.e. of the distance to statistical independence.

The BGD and CL processes are special examples of a general class of pair-interaction processes. The paper will study these processes in the space-homogeneous case. Having in mind the special examples of the BDG and CL processes, we assume that the particle velocities are two-dimensional vectors of constant norm. However, this assumption could be easily waived. The main theorem is that the chaos propagation property is true for these pair-interaction processes. 
The derivation of the BDG and CL kinetic equations follow from this theorem. In the case of the BDG operator, we recover the collision operator or \cite{BDG}. The proof of the theorem generalizes some of the combinatorial arguments of Mark Kac \cite{Kac}. 

This result can be seen as paradoxical at first sight. Indeed, the BDG and CL processes build-up correlations, in the sense that particles tend to eventually become close to each-other (in velocity space). This correlation build-up must be present in any model that is to display ``swarming'' or ``flocking'' behavior, and it  should lead to a breakdown of the statistical independence of the particles. The resolution of this apparent paradox is seen  through an inspection of the time scales. Indeed, the theorem is only valid on any finite time interval at the kinetic scale. At this scale, the average number of collisions undergone by each particle on a finite time interval is finite and uniformly bounded. If $N$ is large and one selects a typical pair of particles and a finite time $t$, none of the particles that have collided with the first particles by that time have also collided  with any particle that has collided with the other particle at time $t$. If one fixes $N$, and lets $t$ become large, pair correlations do develop. However, if one fixes $t$ and lets $N$
 become large, the amount of pair correlation built up by time $t$ goes to zero as $N$ increases.

At equilibrium, the solution of the master equation converges to the so-called invariant density. Due to the special properties of the CL process, it is possible to provide closed expressions for the marginals of the invariant density. When the noise distribution is appropriately scaled as $N \to \infty$, we show that the invariant density cannot be chaotic. Indeed, while the single particle marginal density is uniform, the two-particle marginal density is not. Therefore, the two-particle marginal density is not a tensor product of two copies of the single-particle marginal density, as it should if it would be chaotic. Numerical simulations, reported in a forthcoming work \cite{CCDW}, confirm this result. 
Again, it may seem paradoxical that the invariant measure has strong pair correlations and yet chaos is propagated, and again, it is consideration of the time scales that 
resolves the paradox:  The time for the $N$ particle system to approach equilibrium grows so rapidly with $N$ that the stationary measure in not relevant to time evolution
on the kinetic scale. 

The message of the paper is that the chaos property may be true even for processes that seemingly build-up correlations. However, 
in this case, the correlation build-up capacities of the processes under considerations only  manifest themselves at scales which are large compared to the kinetic scale. To describe these systems at these large scales, kinetic theory is not valid anymore, and alternate theories must be devised. So far, the subject is widely open in the literature and constitute a fascinating area of research.  

We conclude this section by a few more bibliographical remarks. In \cite{DM}, an alternate kinetic model for the Vicsek system has been proposed. It consists of a nonlinear Fokker-Planck equation. It has been derived from a mean-field limit of the Vicsek system in \cite{BCC2}. In biological swarm modeling, most of the authors make use of particle (aka 'Individual-Based') models (see e.g. \cite{Aoki, Chuang, Couzin, Cucker_Smale, Vicsek}) and sometimes, fluid-like hydrodynamic models (see e.g. \cite{Chuang, Mogilner1, TB2}). The use of kinetic models is more rare. There is a kinetic version of the Cucker-Smale model \cite{CFRT}. This kinetic Cucker-Smale model takes the form of a nonlinear non-local Fokker-Planck equation which has been rigorously derived from the mean-field limit of the discrete Cucker-Smale model in \cite{BCC1}. Kinetic models have also been proposed in the context of fish schools \cite{DM2}, bacteria and cell motion \cite{hillen_diffusion_2000, painter_modelling_2009} and ant-trail formation \cite{BDM}. In most cases, their justification is purely formal. The present paper is the first step towards a justification of kinetic models in biological swarm modeling.  

The paper is organized as follows. In section \ref{sec_motivation}, we derive the master equations for the two examples of biological swarm models that we will consider, the BDG and the CL processes. We also provide the definition of the most general pair interaction driven master equation. Section \ref{sec_propa_chaos} is devoted to the proof of the chaos propagation theorem for pair interaction driven master equations and its application to the BDG and CL dynamics.  This proof is modeled in part on the orginal approach of Kac, but with some  differences.  Section \ref{sec_invariant} investigates the invariant measure for the CL process and shows that, under some suitable scaling of the noise distribution, it violates the chaos property. A conclusion is drawn in section \ref{sec_conclu}.

\setcounter{equation}{0}
\section{Motivation: biological swarm models } 
\label{sec_motivation}

\setcounter{equation}{0}
\subsection{Two examples} 
\label{sec_examples}

We consider a population of $N$ {\em agents}, which to be concrete, we take to be fish in a shallow pond, each swimming  at unit speed. In
this work, we are only concerned with the evolution of the velocities, and we neglect the precise spatial location of the fish;
that is, we assume that all $N$ under consideration are sufficiently close to interact with one another. 

The shallow pond is essentially a planar domain, and so the individual velocity vectors belong to
the unit circle ${\mathbb S}^1$. The state space of the system is therefore the torus
$${\mathbb T}_N = [{\mathbb S}^1]^{\times N}\ .$$
The state of the swarm, or school, is then specified by giving a vector
$$\vec v = (v_1,\dots,v_N) \in {\mathbb T}_N\ .$$
All $N$ agents, or fish, are considered as part of a local cluster, and all are interacting with one another.
The evolution of  $\vec v =(v_1,\dots,v_N)$ will be modeled in various ways, all based on the following general scheme: there is a steady Poisson stream of jump times, at which a pair $(i,j)$ is selected at random from $\{1,\dots,N\}$, and then 
these two fish adjust their velocities in some way
\begin{equation}\label{update}
(v_i,v_j) \longrightarrow  (v_i',v_j')\ .
\end{equation}
To complete the specification of the dynamics, we need to give the precise rule for updating the velocities in (\ref{update}).
Here are the rules we shall consider.

\medskip

\noindent{\bf (1) BDG dynamics:}   This rule is designed to lead to a kinetic model first investigated by Bertin, Droz and Gregoire \cite{BDG}.
The idea is that the pair of agents adjust their velocities cooperatively to achieve the same direction of motion, apart from some noise
in their adjustments. More precisely,
define
\begin{equation}\label{bar}
\overline{v}_{i,j} = \frac{v_i+v_j}{|v_i+v_j|}\ .
\end{equation}
Now think of  ${v}_{i,j} \in {\mathbb S}^1$ as unit complex number. Let $W_i$ and $W_j$ be two more unit complex numbers, chosen independently
at random
from a probability distribution $g(w){\rm  d}w$ on ${\mathbb S}^1$, and define
\begin{equation}
\label{updateBDGM}
v_i' = W_i\overline{v}_{i,j}  \qquad{\rm and}\qquad v_j' = W_j\overline{v}_{i,j} \ .
\end{equation}
Regarding ${\mathbb S}^1$ as the unit circle in the complex plane, $W_j\overline{v}_{i,j} $ simply means the product in the complex plane of the random variables $W_i$ and $\overline{v}_{i,j} $, and likewise for the other term. 
Thus, if we write
$$W_i= e^{i\Theta} \quad {\rm and}\quad \overline{v}_{i,j} = e^{i\theta},\quad {\rm  then}\quad W_i\overline{v}_{i,j} = e^{i(\Theta+\theta)}$$
so that the noise is additive in the angles. In the case of no noise,
$W_i$ and $W_j$ are the constant random variable $W_i = W_j =1$, and in the case of small noise, they are random variables that
are strongly peaked around $1$. We suppose that $g(w)$ is symmetric; i.e., 
\begin{equation}
g(w) = g(w^*), 
\label{eq_g_symmetric}
\end{equation}
where $w^*$ denote the complex conjugate of $w$, and that $g(w)$ is somewhat peaked near $w =1$.  

\begin{remark}
This rule could be referred to as the ``Maxwellian BDG'' dynamics, in reference to the fact that the selection of the pair $(i,j)$ is independent of their relative velocity, like the Maxwellian molecular interaction in rarefied gas dynamics. A more general setting would make the collision probability of the pair $(i,j)$ depend on their relative velocity $v_i v_j^*$, but this will be discarded here, for reasons developed below (see remark \ref{rem_pair}). 
\label{rem_Maxwellian_BDG}
\end{remark}

\medskip

\noindent{\bf (2) ``Choose the leader'' (CL) dynamics:} In this variant, 
one of the two agents in the pair decides to adopt the other agent's velocity, though it does not get this velocity exactly right: The new velocity it adopts is the velocity of the other agent up to some noise term. The process is written as follows: if the pair selected is $(1,2)$, and agent 1 decides to adopt the velocity of agent 2, then the velocity
of agent 1 is updated as follows:
$$v_1 \to v_1' := Zv_2\ ,$$
where $Z$ is an independent random variable with values in ${\mathbb S}^1$ and probability $g$. 
As before, we regard ${\mathbb S}^1$ as the unit circle in the complex plane, and $Zv_1$ simply means the product in the complex plane of the random variables $Z$ and $v_1$. Here again, we assume that $g$ is symmetric and satisfies (\ref{eq_g_symmetric}) for the sake of simplicity. We also have in mind that $g$ is peaked around $1$.

We use a fair coin toss, modeled by a Bernoulli variable $B$,  to decide which agent adopts the velocity of the other, and then we have the following description of the jump when pair $(i,j)$
The selected pair of velocities is updated according to
\begin{equation}\label{jump}
\begin{array}{lcl}
v_i' &=& B v_i + (1 - B) Z v_j \\
v_j' &=& B Z v_i + (1-B) v_j
\end{array}
\end{equation}
and all other velocities are unchanged. 

\bigskip
There are many other variants on these basic examples, but for now, let us focus on these two and seek a passage from this
 description of the interactions of individual agents to an evolution equation for the statistical distribution of the velocities in the system.
 For this we shall employ methods of kinetic theory that have been developed for a similar problem concerning colliding molecules in a gas.
 We shall use a probabilistic approach of Marc Kac \cite{Kac}, using a so-called {\em Master equation}.

\setcounter{equation}{0}
\subsection{Master Equations}
\label{subsec_derivation}

\subsubsection{General framework}
\label{subsubsec_general}

We now derive  master equations describing the evolution of the probability density for the state of the system as it undergoes our stochastic processes. 
An advantage with starting at the microscopic level, i.e., the level of individual agents, is that the modeling is much clearer before any large $N$ limits are taken. 

We note that both the BDG or the CL dynamics are Markovian. In general, we consider a Markov process on ${\mathbb T}_N$ and we denote by $\vec V_k \in {\mathbb T}_N$ its state just after the $k$th jump. We define its Markov transition operator $Q$ as usual by
$$Q\varphi(\vec v) = {\rm E}\{\varphi(\vec V_{k+1})\ |\ \vec V_k = \vec v\ \}\ ,$$
for any continuous test function $\varphi$ on ${\mathbb T}_N$. 

Now let $F_k(\vec v)$ denote the probability density of $\vec V_k$ (with respect to the uniform measure on ${\mathbb T}_N$). Then by definition, one has
$${\rm E}(\varphi(\vec V_{k+1})) = \int_{{\mathbb T}_N}\varphi(\vec v)F_{k+1}(\vec v){\rm d}^Nv\ .$$
On the other hand, by standard properties of the conditional expectation,
$${\rm E}(\varphi(\vec V_{k+1}))  = {\rm E} ({\rm E}\{ \varphi(\vec V_{k+1})\ |\ V_k \} ) = 
 \int_{{\mathbb T}_N}Q\varphi(\vec v)F_{k}(\vec v){\rm d}^Nv\ .$$
That is,
 \begin{equation}\label{adj}
 F_{k+1} = Q^*F_k, 
 \end{equation}
where $Q^*$ is the adjoint of $Q$ in $L^2({\mathbb T}_N,{\rm d}^Nv)$. 

The next step is to construct a time continuous process which will lead to a time-continuous master equation. The state of the process is now a function of time $\vec v (t) \in {\mathbb T}_N$ and the probability density $F(\vec v,t)$ is a function of the continuous time parameter $t$ instead of the discrete jump index $k$. We assume that $F(\vec v,t+dt)$  only depends on $F(\vec v,t)$ and not on the past values $F(\vec v,s)$ for  $s<t$. In this way, we can construct a time-continuous
Markov process. Thus, we assume that between time $t$ and $t + dt$ the probability that a
collision occurs is
$$ \lambda Q^* F (\vec v,t) \, dt + o(dt), $$
where $\lambda$ is a constant. We assume that the probability that multiple collisions occur in the
time interval $[t, t+dt]$ is negligible. Thus, since there are $N$ particles colliding independently, we have, for $dt$: 
$$ F(\vec v, t + dt) = N \lambda Q^* F(\vec v, t) \, dt + (1 -  N \lambda F(\vec v, t)) \, dt + o(dt) .$$
Equivalently, we can write 
$$ F(\vec v, t + dt) - F(\vec v, t) = N \lambda [ Q^* F( \vec v, t) \, dt - F(\vec v, t)] \, dt + o(dt) , $$
which leads to the Master Equation in the limit $dt \to 0$: 
$$ \frac{d}{dt} F(\vec v, t) = \lambda N [ Q^* - I] F(\vec v, t), $$
where $I$ is the identity operator. This equation must be complemented by the initial condition
$$ F(\vec v, 0) = F_0(\vec v) $$
where $F_0$ is the initial probability distribution. In the remainder, we scale time in such a way that $\lambda = 1$ and we define
\begin{equation}
 L = N [ Q - I] \quad \mbox{ and } \quad L^* = N [ Q^* - I]. 
\label{eq_def_LL*_gene}
\end{equation}
We summarize the previous discussion in the following definition: 

\begin{definition} The time-continuous master equation associated to a discrete Markov process of transition operator $Q$ is written  
\begin{eqnarray}
&& \hspace{-1cm}
\frac{d}{dt} F(\vec v, t) = L^* F(\vec v, t), 
\label{eq_master_gene} \\
&& \hspace{-1cm}
F(\vec v, 0) = F_0(\vec v) . \label{eq_master_gene_ini}
\end{eqnarray}
where $L^* = N [ Q^* - I]$ with $Q^*$ the adjoint of $Q$ and $F_0$ is the initial probability distribution. 
\label{def_master}
\end{definition}

\noindent
In the next sections, we determine the master equations of the BDG and CL processes successively.

\subsubsection{The BDG dynamics}
\label{subsubsec_master_BDG}

We state the following:

\begin{proposition}
The master equation of the BDG dynamics is written (\ref{eq_master_gene}) with 
\begin{equation}
L^* = N \left(\begin{matrix} N\cr 2\cr \end{matrix}\right)^{-1}\sum_{i<j}(Q_{(i,j)}^* - I)\ ,
\label{eq_def_L*_BDG}
\end{equation}
and the binary interaction operator $Q_{(i,j)}^*$ given by: 
\begin{equation}
Q^*_{(i,j)} F(\vec v) = 
\int_{{\mathbb T}^2} F(v_1,\dots,y_i,\dots,y_j,\dots,v_N) \, g\left(v_i\overline{y}_{i,j}^*\right) \, g\left(v_j\overline{y}_{i,j}^*\right) \, {\rm d}y_i{\rm d}y_j\ .
\label{eq_def_Qij*_BDG}
\end{equation}
\label{prop_master_BDG}
\end{proposition}

\medskip
\noindent
{\bf Proof:}
>From the definition (\ref{updateBDGM}) of the BDG dynamics, we get
\begin{eqnarray*} 
Q\varphi (\vec v) & = &  {\rm E}\{\varphi(\vec V_{k+1})\ |\ \vec V_k = \vec v\ \} \\
& = &
{\rm E} \varphi(v_1,\dots, W_i\overline{v}_{i,j}, \dots, W_j\overline{v}_{i,j},\dots,v_N) \\
& = & \frac{2}{N(N-1)} \sum_{i<j} \int_{{\mathbb T}_2} \varphi(v_1,\dots, w_i\overline{v}_{i,j},\dots, w_j\overline{v}_{i,j},\dots, v_N) \, g(w_i) \, g(w_j) \, {\rm d}w_i
{\rm d}w_j \ .
\end{eqnarray*}
To compute the adjoint of $Q$, we note that for any probability density $F$ on ${\mathbb T}_N$,
\begin{eqnarray*}
&&\int_{{\mathbb T}_N}F(v_1,\dots,v_N) \times \\
&&\left[\int_{{\mathbb T}_2}\varphi(v_1,\dots, w_i\overline{v}_{i,j},\dots, w_j\overline{v}_{i,j},\dots, v_N) \, g(w_i) \, g(w_j) \, {\rm d}w_i
{\rm d}w_j \right] {\rm d}v_1\dots {\rm d}v_N =\\
&& = \int_{{\mathbb T}_N} \left[\int_{{\mathbb T}_2} F(v_1,\dots,v_i,\dots,v_j,\dots,v_N)  \, g\left(y_i\overline{v}_{i,j}^*\right) \, g\left(y_j\overline{v}_{i,j}^*\right) \, {\rm d}v_i{\rm d}v_j\right] \times \\
&& \hspace{3cm} 
\varphi(v_1,\dots , y_i,\dots, y_j,\dots, v_N) \, {\rm d}v_1\dots  {\rm d}y_i, \dots  {\rm d}y_j\dots, {\rm d}v_N, 
\end{eqnarray*}
where we have introduced the variables
$$y_i =  w_i\overline{v}_{i,j} \qquad {\rm and} \qquad y_j = w_j\overline{v}_{i,j}\ .$$
Now changing the names of variables, we finally have
\begin{multline}
Q^*F(v_1,\dots,v_N) =\\ \frac{2}{N(N-1)}\sum_{i<j}
\int_{{\mathbb T}_2} F(v_1,\dots,y_i,\dots,y_j,\dots,v_N)g\left(v_i\overline{y}_{i,j}^*\right)g\left(v_j\overline{y}_{i,j}^*\right){\rm d}y_i{\rm d}y_j\ ,
\end{multline}
where
$$\overline{y}_{i,j} = \frac{y_i+y_j}{|y_i+y_j|}\ .$$
Therefore, we can write
$$ Q^* = \left(\begin{matrix} N\cr 2\cr \end{matrix}\right)^{-1} \sum_{i<j} Q_{(i,j)}^*, $$
with $Q_{(i,j)}^*$ defined by (\ref{eq_def_Qij*_BDG}). Then, using (\ref{eq_def_LL*_gene}), we get eq. (\ref{eq_def_L*_BDG}). \endproof

\begin{remark}
It is useful to note that the adjoint $L$ of $L^*$ which corresponds to the Markov transition operator, is defined by
\begin{equation}
L = N \left(\begin{matrix} N\cr 2\cr \end{matrix}\right)^{-1} \sum_{i<j}(Q_{(i,j)} - I)\ ,
\label{eq_def_L}
\end{equation}
with 
$$
Q_{(i,j)}\varphi(\vec v) = \\ \int_{{\mathbb T}_2} \varphi(v_1,\dots, w_i\overline{v}_{i,j},\dots, w_j\overline{v}_{i,j},\dots, v_N) \, g(w_i) \, g(w_j) \, {\rm d}w_i {\rm d}w_j\ .$$
\end{remark}

 \subsubsection{The CL dynamics}
 \label{subsubsec_master_CL}
 
We now derive the master equation for the CL model. We introduce the notation  $(v_1,\dots,\widehat{v_i},\dots,v_N)$ for the $n-1$ tuple formed by removing $v_i$ from $\vec v$ and
$$[F]_{\widehat{i}} (v_1,\dots,\widehat{v_i},\dots,v_N) :=  \int_{{\mathbb S}^1} F(v_1,\dots,v_N){\rm d}v_i$$
for the marginal of $F$ obtained by integrating $v_i$ out. 
We show the
 
\begin{proposition}
The master equation of the CL dynamics is written (\ref{eq_master_gene}) with $L^*$ given by (\ref{eq_def_L*_BDG}) and the binary interaction operator $Q_{(i,j)}^*$ by: 
\begin{equation}
Q_{(i,j)}^*F(\vec v) = \frac{1}{2}\left[
[F]_{\widehat{i}}(v_1,\dots,\widehat{v_i},\dots,v_N) +  [F]_{\widehat{j}}(v_1,\dots,\widehat{v_j},\dots,v_N) \right] \, g(v_i^*v_j) \ .
\label{eq_def_Qij*_CL}
\end{equation}
\label{prop_master_CL}
\end{proposition}

\noindent
{\bf Proof:} From the definition (\ref{jump}) of the CL process, the Markov transition operator $Q$ is given by:
\begin{eqnarray}
&&\hspace{-1cm}
Q\varphi(\vec v) = \frac{1}{N(N-1)} \sum_{i<j} \int_{{\mathbb S}^1}\left[
\varphi(v_1,\dots, zv_j,\dots, v_j,\dots, v_N) \right. \nonumber \\
&&\hspace{4cm}
\left. + 
\varphi(v_1,\dots, v_i,\dots, zv_i,\dots, v_N) \right] \, g(z) \, {\rm d} z\ .
\label{Qdef}
\end{eqnarray}
To compute the adjoint of $Q$, we note that for any probability density $F$ on ${\mathbb T}_N$,
\begin{multline}\label{cal1}
\int_{{\mathbb T}_N}F(v_1,\dots,v_N)\left[\int_{{\mathbb S}^1}\varphi(v_1,\dots, zv_j,\dots, v_j,\dots, v_N)g(z){\rm d}z\right]{\rm d}v_1\dots {\rm d}v_N = \\
= \int_{{\mathbb T}_{N-1}} [F]_{\widehat{i}} (v_1,\dots,\widehat{v_i},\dots,v_N)
\int_{{\mathbb S}^1}
\varphi(v_1,\dots , zv_j,\dots, v_j,\dots, v_N) \, g(z)\\
\hspace{4cm} {\rm d}z \, {\rm d}v_1\dots \widehat {{\rm d}v_i}, \dots, {\rm d}v_N
\end{multline}
We next introduce a new variable 
$$y_i = zv_j \quad \mbox{ (or equivalently }  z = v_j^*y_i \mbox{ )}. $$
Evidently ${\rm d}z = {\rm d}y_i$. Additionally, since $v_i$ has been integrated out, it has disappeared from (\ref{cal1}). Therefore, we can change the name $y_i$ into $v_i$ without any confusion. We can then rewrite (\ref{cal1}) as
\begin{multline}\label{cal2}
\int_{{\mathbb T}_N} F(v_1,\dots,v_N)
\left[\int_{{\mathbb S}^1}\varphi(v_1,\dots, zv_j,\dots, v_j,\dots, v_N)g(z){\rm d}z\right]
{\rm d}v_1\dots {\rm d}v_N = \\
= \int_{{\mathbb T}_N}[F]_{\widehat{i}}(v_1,\dots,\widehat{v_i},\dots,v_N) \, g(v_j^*v_i) \, 
\varphi(v_1,\dots , v_i,\dots, v_j,\dots, v_N) \\
\hspace{4 cm} {\rm d}v_1\dots {\rm d}v_i, \dots, {\rm d}v_N\ .
\end{multline}
Using this formula and its analog for $i$ and $j$ exchanged, we see that the master equation for the CL dynamics is given by (\ref{eq_def_LL*_gene}), (\ref{eq_def_L*_BDG}) with $Q_{(i,j)}$ given by (\ref{eq_def_Qij*_CL}). \endproof

\begin{remark}
Again, we note that the adjoint $L$ of $L^*$ is defined by (\ref{eq_def_L}) with $Q_{(i,j)}$, the adjoint of $Q_{(i,j)}^*$, given by:
$$
Q_{(i,j)}\varphi(\vec v)= \int_{{\mathbb S}^1}\varphi(v_1,\dots, zv_j,\dots, v_j,\dots, v_N)g(z){\rm d}z\ .
$$
\end{remark}

\subsection{Extension: Pair-Interaction driven Master Equation}
\label{subsec_master_PI}
 
The master equations of the BDG and CL dynamics are two examples of a class of master equations which we will call 'Pair Interaction driven Master Equations', defined below. 
 
\begin{definition}[Pair Interaction Driven Master Equation] A {\em pair interaction driven Master equation} is an equation of the form
$$\frac{\partial}{\partial t}F(\vec v,t) = L^*F(\vec v,t),$$
describing the evolution of probability densities  on some product space $X_N$ with elements $\vec v = (v_1,\dots,v_N)$
where
 $$L^* = N \sum_{i<j} p_{i,j}(\vec v) \, (Q_{(i,j)}^* - I)\ .$$
The operators  $Q_{(i,j)}$ are Markov operators on functions on $X_N$ such that  $Q_{(i,j)}\varphi = \varphi$ whenever $\varphi$ does not depend on 
either $v_i$ or $v_j$. The pair selection probabilities $p_{i,j}(\vec v)$ are such that $p_{i,j}(\vec v) \geq 0$ and
$$\sum_{i<j}p_{i,j}(\vec v) =1\ .$$
\label{def_master_pair}
\end{definition}

We have given two examples already: the CL and BDG master equations. The {\em Kac Master equation} \cite{Kac} is another example and is described below.
 
\begin{example} The Kac Master equation. In this example, $X_N$ is the sphere in ${\mathord{\mathbb R}}^N$ of radius $\sqrt{N}$, 
\begin{eqnarray*}
\hspace{-1cm} & & 
Q_{(i,j)}\varphi ({\vec v}) =  \int_{-\pi}^{\pi}\rho(\theta)
\varphi(R_{i,j,\theta}\vec v){\rm d}\theta\ ,\\
\hspace{-1cm} & & 
R_{i,j,\theta}\vec v = (v_1,v_2,\dots,\cos\theta \, v_i + \sin\theta \, v_j,\dots,-\sin\theta \, v_i+\cos\theta \, v_j,\dots,v_n),  \\
\hspace{-1cm} & & 
p_{i,j} = \frac{2}{N (N-1)}
\end{eqnarray*}
and $\rho$ is a probability density on ${\mathbb S}^1$. 
The operators $Q_{(i,j)}$ in the Kac model are self adjoint with respect to the uniform probability measure on the sphere $S^{N-1}$, which is therefore the invariant measure for this process. In other words the Kac process is {\em reversible} meaning that is satisfies {\em detailed balance}: if you saw a movie of the process running backwards, there would be no clue that it was running backwards. 
\end{example}

By contrast to the Kac model, the BDG and CL models {\em do not} have detailed balance property and time reversibility.  If you ran the movie backwards, you would see pairs of fish with similar velocities changing them to differ in a random way.  For these non-reversible processes, it is not so easy to determine the invariant measure, though it will exist and be unique for each $N$ for our processes under mild assumptions on the noise distribution $g$. In section \ref{sec_invariant}, it will be possible to determine the marginals of the invariant density of the CL dynamics in closed form. However, this simplification is not possible for the BDG dynamics. 

In the next section, we show that pair interaction driven master equations with uniform selection probabilities $p_{i,j} = 2/(N(N-1))$ do have the propagation of chaos property, and therefore, satisfy a kinetic equation at the kinetic time scale. However, in section \ref{sec_invariant}, we show that the equilibrium density of the CL dynamics cannot satisfy the propagation of chaos property, meaning that this property may break down at larger time scales.

\setcounter{equation}{0}  
\section{Propagation of Chaos}
\label{sec_propa_chaos}

\subsection{Definition}
\label{subsec_CP_definition}

To pass to a kinetic description, and then on to a hydrodynamic description, the key step is a propagation of chaos result. That may seem unlikely in the cases of the CL and BDG dynamics which are expected to build pair correlations. However, the time scales at which pair correlations built up may be longer than the kinetic time scale at which a kinetic model is expected to be valid. In the present section, we shall see that chaos is propagated in both the BDG and CL models at the kinetic time scale. In section \ref{sec_invariant}, we prove that the invariant measure of the CL dynamics is not chaotic; it exhibits pair  correlation. These two observations are not self-contradictory since propagation of chaos holds only on a finite time scale while the invariant measure is reached as time tends to infinity.

\begin{definition}[Chaos]   Let $X_N$ be the $N$-fold cartesian product of a polish space $X$ equipped with some reference measure $\mu$. 
Let $f$ be a given probability density on $X$.  For each $N\in {\mathord{\mathbb N}}$, let $F_{(N)}$ be a probability density on $X_N$ with respect to $\mu^{\otimes N}$. 
The sequence $\{F_{(N)}\}_{N\in {\mathord{\mathbb N}}}$ of probability densities on $X_N$ is {\em $f$-chaotic}
in case
 
\medskip
\noindent{\it(1.)}  Each $F_{(N)}$  is a symmetric function of  
$\{v_1,v_2,\dots,v_N\}$ 

\medskip
\noindent{\it(2.)} For each fixed $k$, and any bounded measurable function $\phi$ on ${\mathord{\mathbb R}}^k$,
\begin{eqnarray*}
&&\hspace{-1cm} 
\lim_{N\to\infty}\int_{X_N}\phi(v_1,v_2,\dots,v_k)F_{(N)}(\vec v){\rm d}\sigma =  \int_{X^k}\phi(v_1,v_2,\dots,v_k)\prod_{j=1}^k f(v_j){\rm d}^kv\ .
\end{eqnarray*}
\label{def_CP}
\end{definition}

\medskip
\noindent
Kac \cite{Kac} proved that the semigroup $e^{tL^*}$ associated to Kac's master equation {\em propagates chaos}. More precisely, {\bf Kac's Theorem} is stated as follows: 

\begin{theorem}[Propagation of chaos]
Let $\{F_{(N)}\}_{N \in {\mathbb N}}$ be $f$--chaotic. Then the family $\{e^{tN(Q-I)} F_{(N)}\}_{N \in {\mathbb N}}$ is
$f(\cdot,t)$--chaotic where
 $f(v,t)$ is the solution of
\begin{equation}\label{kke}
\frac{\partial f}{\partial t}(v,t) = Q(f,f)(v,t)\qquad{\rm with}\qquad f(v,0) = f(v)\ ,
\end{equation}
with
$$
Q(f,f) = 2\int_{{\mathord{\mathbb R}}}\int_{-\pi}^{\pi} \left[
f(v')f(w') - f(v)f(w)\right] \, \rho(\theta) \, {\rm d} \theta {\rm d} w \ ,
$$
and 
$$ v' = \cos\theta \, v + \sin\theta \, w, \quad w' = -\sin\theta \, v +\cos\theta \, w. $$
\end{theorem}

\noindent
Eq. (\ref{kke}) is called the {\em Kac-Boltzmann equation}.  In this section we prove a propagation of chaos result valid in the general class of pair interaction driven Master equations. We shall use this result to discuss the kinetic limits of the BDG and CL dynamics.

\subsection{Propagation of chaos for pair interaction driven master equations}
\label{subsec_CP_PIM}

Consider a general Master equation
\begin{equation}
\frac{\partial}{\partial t}F = L^*F,
\label{eq_GME}
\end{equation}
for a probability density $F$ on ${\mathbb T}_N$ of the form
\begin{eqnarray}
&&\hspace{-1cm}
L^*F = N (Q^* - I) F = \frac{2}{N-1}\sum_{i<j}(Q^*_{(i,j)} -I)F, 
\label{eq_GME_form} 
\end{eqnarray}
where $Q_{i,j}$ is a Markovian operator acting on $F$ through $v_i$ and $v_j$ alone. The goal of this section is to prove the following: 

\begin{theorem} Let $\{F_{0,N}^{(N)}\}_{N \in {\mathbb N}}$ be $f$-chaotic.  Then for each $t>0$, the family of marginals 
$\{e^{tL^*}F_{0,N}^{(N)}\}_{N \in {\mathbb N}}$ associated to eq. (\ref{eq_GME}), where $L^*$ is a pair-interaction operator of the form (\ref{eq_GME_form}), is $f(\cdot,t)$-chaotic where $f(\cdot,t)$ satisfies the following Boltzmann equation: 
\begin{equation}
\frac{{\partial}}{{\partial t}}f(v,t) = 2\left[\int_{{\mathbb S}^1} Q^*_{(1,2)}f^{\otimes 2}(v,w) \, {\rm d} w -  f(v,t)\right]\ ,
\label{eq_Boltzmann_PIM}
\end{equation}
associated to the initial condition $f(v,0) = f$. 
\label{thm_chaos_prop}
\end{theorem}

\noindent
Before proving Theorem \ref{thm_chaos_prop}, we make some preliminary comments. 
Let initial data $F_{0,N}$ be given, and let us compute the evolution of $F^{(1)}_t$, the single particle marginal at time $t$. For any test function $\varphi(v)$ of the single coordinate 
$v \in {\mathbb S}^1$, we have
\begin{eqnarray*}
\int_{{\mathbb S}^1} \varphi(v) \, F^{(1)}_t(v) \, {\rm d} v &=& \int_{{\mathbb T}_N} \varphi(v_1) \, e^{tL^*}F_{0,N}(\vec v) \, {\rm d} v_1 \ldots{\rm d} v_N \\
&=& \int_{{\mathbb T}_N} e^{tL}\varphi(v_1) \, F_{0,N}(\vec v) \, {\rm d} v_1 \ldots {\rm d} v_N\ . 
\end{eqnarray*}
A similar relation holds for the two particle marginal and so on. So, to study the evolution of low dimensional marginals, it is helpful to understand the behavior of expressions of the form $e^{tL}\varphi(\vec v)$ when $\varphi(\vec v)$ depends on only finitely many coordinates in $\vec v$. 
 
It is clear that in general, for a bounded continuous function $\varphi$ on ${\mathbb T}_N$,
$$\|L \varphi\|_\infty \leq 2N\|\varphi\|_\infty\ .$$
However, if $\varphi$ depends only on $v_1,\dots,v_k$, tighter bounds are valid. This is because 
$$i,j > k \quad \Rightarrow \quad Q_{(i,j)}\varphi = \varphi\ ,$$
and so
\begin{equation}
L \varphi = \frac{2}{N-1} \sum_{i<j}  (Q_{(i,j)} -I) \varphi = \frac{2}{N-1} \sum_{i=1}^k \sum_{j=i+1}^N (Q_{(i,j)} -I)\varphi \ ,
\label{eq_Lphi_v1}
\end{equation}
and thus, as soon as $k \geq 1$: 
\begin{equation}
\|L\varphi\|_\infty \leq \frac{2}{N-1} \, k (N-\frac{k+1}{2}) \, 2 \|\varphi\|_\infty \leq 4k\|\varphi\|_\infty \ .
\label{eq_fundamental_property}
\end{equation}
We can now state the following fundamental lemma:
 
\begin{lemma} Let $\varphi$ be a function depending only on $v_1,\dots,v_k$. We can regard $\varphi$ as a function on
$X_N$ for each $N\in {\mathord{\mathbb N}}$, $N \geq k$. Then,
the power series
\begin{equation}
e^{tL}\varphi = \sum_{k=0}^\infty \frac{t^k}{k!}L^k \varphi,
\label{eq_etL}
\end{equation}
converges absolutely in $L^\infty$, uniformly in $N \in {\mathbb N}^*$  and $t \in [0,T]$ for any $T<1/4$. 
\label{lem_absolute_convergence}
\end{lemma}
 
\noindent{\bf Proof:} 
Consider first the case in which $\varphi$ depends only on one variable. Without loss of generality, owing to the permutation symmetry of the problem, we can set this variable to $v_1$. Then from (\ref{eq_Lphi_v1}), $L\varphi$ is an average of functions depending on only two velocities. Likewise,  $L^2\varphi$ is a combination of terms depending only on three velocities and so on.  By what we have noted above, we can expect the following formula:
\begin{equation}
\|L^k\varphi\| \leq  4^{k} \, k! \, \|\varphi\|_\infty\ .
\label{eq_||Lkphi||}
\end{equation}

To show (\ref{eq_||Lkphi||}) we prove that $L^k \varphi$ is of the form
\begin{equation}
L^k \varphi = \left( \frac{2}{N-1} \right)^k  \sum_{s \in S_k} \psi_s^{(k)}, 
\label{eq_Lkphi_induction}
\end{equation}
where the set $S_k$ is a set of multi-indices $s=(1,s_1, \ldots, s_k)$, such that
\begin{equation}
\mbox{Card} \, S_k \leq  k! \, \prod_{j=1}^k (N - \frac{j+1}{2}). 
\label{eq_Card_Sk}
\end{equation}
The function $\psi_s^{(k)}$ depends only on the $k+1$ variables $(v_1, v_{s_1}, \ldots , v_{s_k})$ and satisfies
\begin{equation}
\| \psi_s^{(k)} \|_\infty \leq 2^k \|\varphi\|_\infty.
\label{eq_||psi_s^k||}
\end{equation}
Of course, (\ref{eq_||Lkphi||}) results from (\ref{eq_Lkphi_induction}), (\ref{eq_Card_Sk}) and (\ref{eq_||psi_s^k||}) and the remark that 
$$ \frac{1}{(N-1)^k}  \prod_{j=1}^k (N - \frac{j+1}{2}) \leq 1. $$
We note that, if $k+1>N$, some of the variables $(v_1, v_{s_1}, \ldots , v_{s_k})$ may be identical, but the arguments below are still valid in this case.  

The proof of (\ref{eq_Lkphi_induction}) is by induction. For $k=1$, using (\ref{eq_Lphi_v1}), we have
$$ L \varphi = \frac{2}{N-1} \sum_{j=2}^N (Q_{(1,j)} - I) \varphi. $$
Therefore, letting $S_1 = \{2, \ldots, N\}$ and $\psi_s^{(1)} = \psi_j^{(1)} = (Q_{(1,j)} - I) \varphi$, we can write $L \varphi$ according to formula (\ref{eq_Lkphi_induction}). Clearly, Card $S_1 = N-1$ in accordance to (\ref{eq_Card_Sk}). Finally, by the Markov property, we have $\| \psi_j^{(1)} \|_\infty \leq 2 \| \varphi \|_\infty$, which is consistent with (\ref{eq_||psi_s^k||}). Therefore, (\ref{eq_||Lkphi||}) is proved for $k=1$. 

Now, we assume that (\ref{eq_||Lkphi||}) is true for $k$ and try to deduce it for $k+1$. 
By the induction hypothesis, $\psi_s^{(k)}$ depends only on $k+1$ variables so (\ref{eq_Lphi_v1}) applies and we have 
\begin{eqnarray*} 
L^{k+1} \varphi &=& \left( \frac{2}{N-1} \right)^k  \sum_{s \in S_k} L \psi_s^{(k)} \\
&=& \left( \frac{2}{N-1} \right)^{k+1}  \sum_{s \in S_k} \, \, \sum_{m \in s, \ell >m} (Q_{(m,\ell)} - I) \psi_s^{(k)} . 
\end{eqnarray*}
The expression $m \in s$ means that $m$ takes any of the indices $\{ 1,s_1, \ldots, s_k \}$ present in the multi-index $s$. We know from the computation in formula (\ref{eq_fundamental_property}) that there are $(k+1) (N-\frac{k+2}{2})$ such pairs $(m,\ell)$ for a given $s$. Defining $S_{k+1}$ as the set of the so-constructed multi-indices $\{s' =  (s,\ell) \}$, we can write 
\begin{eqnarray} 
L^{k+1} \varphi =  &=& \left( \frac{2}{N-1} \right)^{k+1}  \sum_{s' \in S_{k+1}} \psi_{s'}^{(k+1)} , \label{eq_Lk+1phi_induction}
\end{eqnarray}
with 
$$ \psi_{s'}^{(k+1)} = (Q_{(m,\ell)} - I) \psi_s^{(k)}. $$
Clearly, $\psi_{s'}^{(k+1)}$ is a function of the $k+2$ variables $(v_1, v_{s_1}, \ldots , v_{s_k}, v_m)$ and we have 
\begin{eqnarray} 
\mbox{Card} \, S_{k+1} &=& (k+1) (N-\frac{k+2}{2}) \, \mbox{Card} \, S_{k} \nonumber \\
&= & (k+1) (N-\frac{k+2}{2}) \, k! \, \prod_{j=1}^k (N - \frac{j+1}{2}) \nonumber\\
&=& (k+1)! \, \prod_{j=1}^{k+1} (N - \frac{j+1}{2}) . \label{eq_Card_Sk+1}
\end{eqnarray}
Finally, we have, using the Markov property 
\begin{equation}
\| \psi_{s'}^{(k+1)} \|_\infty  \leq 2 \| \psi_s^{(k)} \|_\infty \leq 2^{k+1} \|\varphi\|_\infty.  \label{eq_||psi_s^k+1||}
\end{equation}
Now, collecting (\ref{eq_Lk+1phi_induction}), (\ref{eq_Card_Sk+1}), (\ref{eq_||psi_s^k+1||})
shows that the induction hypothesis is valid at rank $k+1$. This proves (\ref{eq_||Lkphi||}).

Using (\ref{eq_||Lkphi||}), we deduce that
$$\left\Vert \frac{t^k}{k!}L^k \varphi \right\Vert_\infty \leq  (4t)^k \, \|\varphi\|_\infty\ ,$$
uniformly in $N$, and so for $t<1/4$, we have the uniform absolute convergence of the series (\ref{eq_etL}).  
 
If now $\varphi$ is a function $p$ variables, $p \geq 2$, formula (\ref{eq_||Lkphi||}) is changed into 
\begin{eqnarray}
\|L^k\varphi\| &\leq&  4^{k} \, k! \, \left( \begin{array}{c} p+k-1 \\ k \end{array} \right) \,   \|\varphi\|_\infty \label{eq_||Lkphi||_pvar_0} \\
&\leq& 4^{k} \, k! \, \frac{(k+p-1)^{p-1}}{(p-1)!} \,   \|\varphi\|_\infty\  \nonumber \\
&\leq& C_p \, 4^{k} \, k! \, (k+1)^{p-1} \,   \|\varphi\|_\infty\ . 
\label{eq_||Lkphi||_pvar}
\end{eqnarray}
The proof follows the same lines as above. The only thing to note is that Card $S_k$ is now changed into
$$ \mbox{Card} \, S_k = \prod_{j=1}^k (p+j-1) (N- \frac{p+j}{2}) , $$
the other expressions remaining identical. Then, we get  
\begin{equation}
\left\Vert \frac{t^k}{k!}L^k \varphi \right\Vert_\infty \leq C_p \,  (4t)^k \, (k+1)^{p-1}  \, \|\varphi\|_\infty\ ,
\label{tkLkphi}
\end{equation}
uniformly in $N$. For any given $p \geq 2$, the right-hand side of (\ref{tkLkphi}) is still the general term of a convergent series for $t<1/4$. Therefore, the series (\ref{eq_etL}) is still absolutely uniformly convergent for $t<1/4$, which ends the proof of Lemma \ref{lem_absolute_convergence}. \endproof

\begin{remark}
>From the last proof, we note that the series (\ref{eq_etL}) is {\bf not} uniformly convergent with respect to $p$. 
\end{remark}

\noindent
We now have the following 

\begin{lemma}
Suppose that $F_{0,N}$ is a symmetric probability density on ${\mathbb T}_N$, and suppose that $\varphi^{(k)}$ depends only on $\, v_1,\dots,v_k$ and is $L^\infty$. Define $\varphi^{(k+1)}$ by
\begin{equation}
\label{induc}
\varphi^{(k+1)}(v_1,\dots,v_{k+1}) =  2 \sum_{i=1}^k(Q_{(i,k+1)} -I)\varphi^{(k)}(v_1,\dots,v_k)
\ .
\end{equation}
Then, if $N \geq k+1$, we have 
\begin{eqnarray}
&&\hspace{-1cm} 
\int_{{\mathbb T}_N} F_{0,N} \, L\varphi^{(k)} \, {\rm d} v_1\dots{\rm d} v_N = \int_{{\mathbb T}_N}
F_{0,N} \, ( \varphi^{(k+1)} + \tilde \varphi^{(k+1)} ) \, {\rm d} v_1\dots{\rm d} v_N ,
\label{eq_int_FN_Lphik}
\end{eqnarray}
where $\tilde \varphi^{(k+1)}$ only depends on $\, v_1,\dots,v_{k+1}$ and is such that
\begin{equation} \| \, \tilde \varphi^{(k+1)} \, \|_\infty \leq 6 \frac{k \, (k-1)}{N-1}  \, \| \varphi^{(k)} \|_{\infty}. 
\label{eq_estim_SN}
\end{equation}
\label{lem_phi_k+1}
\end{lemma}

\medskip   
\noindent
{\bf Proof:}
We compute, using (\ref{eq_Lphi_v1}) and that $N \geq k+1$: 
\begin{eqnarray*}
&&\hspace{-1cm}
\int_{{\mathbb T}_N} F_{0,N} \, L\varphi^{(k)} \, {\rm d} v_1\dots{\rm d} v_N =  \\
&&\hspace{0.5cm}
= 
 \frac{2}{N-1} \sum_{i=1}^k \sum_{j=i+1}^N \int_{{\mathbb T}_N} F_{0,N} \, (Q_{(i,j)} -I) \varphi^{(k)} \, {\rm d} v_1\dots{\rm d} v_N \\
&&\hspace{0.5cm}
= 
2 \frac{N-k}{N-1} \sum_{i=1}^k \int_{{\mathbb T}_N} F_{0,N} \, (Q_{(i,k+1)} -I) \varphi^{(k)} \, {\rm d} v_1\dots{\rm d} v_N  + \\
&&\hspace{4.5cm}
+ \frac{2}{N-1} \sum_{i<j\leq k}  \int_{{\mathbb T}_N} F_{0,N} \, (Q_{(i,j)} -I) \varphi^{(k)} \, {\rm d} v_1\dots{\rm d} v_N \\ 
&&\hspace{0.5cm}
=  \sum_{i=1}^k \int_{{\mathbb T}_N} F_{0,N} \, ( \varphi^{(k+1)} + \tilde \varphi^{(k+1)} ) \, {\rm d} v_1\dots{\rm d} v_N  ,
\end{eqnarray*}
where we have used the symmetry in the second equality, and where 
\begin{eqnarray}
&&\hspace{-0.5cm}
\tilde \varphi^{(k+1)} (v_1,\dots,v_{k+1}) = 
2 \frac{k-1}{N-1} \sum_{i=1}^k (Q_{(i,k+1)} -I) \varphi^{(k)} + \nonumber\\
&&\hspace{6.5cm}
+ \frac{2}{N-1} \sum_{i<j\leq k} (Q_{(i,j)} -I) \varphi^{(k)} . \label{eq_SN}
\end{eqnarray}
This shows (\ref{eq_int_FN_Lphik}). Now, from the Markov property of $Q_{(i,j)}$ and from (\ref{eq_SN}), we get (\ref{eq_estim_SN}), which ends the proof of the Lemma.  
\endproof

\medskip
\noindent
Now consider any $\varphi^{(m)}$ depending only on $v_1,\dots,v_m$. Since if $F_{0,N}$ is symmetric, so is each $(L^*)^k F_{0,N}$. Therefore, we can repeatedly apply the previous lemma and so on. Using (\ref{induc}), we inductively define $\varphi^{(m+1)}$, $\varphi^{(m+2)}$, \ldots $\varphi^{(m+k)}$. We note that, if $\varphi^{(m)}$ does not depend on $N$, neither does $\varphi^{(m+k)}$. Therefore, the functions $\varphi^{(m+k)}$ are good candidates to express what happens in the limit $N \to \infty$. However, for a given $N$, formula (\ref{eq_int_FN_Lphik}) is only valid until $m+k \leq N$. A special treatment is required for indices $k$ such that $m+k >N$. The contribution of the remainder $\tilde \varphi^{(k+1)}$ needs also to be estimated. These are the goals of the following Lemma.   

\begin{lemma}
We assume that $\{F_{0,N}\}$ is $f$-chaotic (see definition \ref{def_CP}). Then, for $t<1/4$, we have: 
\begin{eqnarray}
&& \hspace{-1cm}
\lim_{N\to\infty} \int_{{\mathbb T}_N} F_{0,N} \, \, e^{tL} \varphi^{(m)} \, {\rm d} v_1\dots{\rm d} v_N = \nonumber \\
&& \hspace{1cm}
= \sum_{k=0}^\infty \frac{t^k}{k!} \int_{{\mathbb T}_{k+m}} \left(\prod_{j=1}^{k+m}f(v_j) \right) \varphi^{(m+k)}(v_1,\dots,v_{k+m}) \, {\rm d} v_1\dots{\rm d} v_{k+m}\ . 
\label{eq_lim_int_FN_etLphi}
\end{eqnarray}
\label{lem_lim_int_FN_etLphi}
\end{lemma}

\noindent
{\bf Proof:} 
>From Lemma \ref{lem_absolute_convergence}, we have
$$ \lim_{N\to\infty} \int_{{\mathbb T}_N}F_{0,N} \, \,  e^{tL}\varphi^{(m)} \, {\rm d} v_1\dots{\rm d} v_N = \lim_{H\to\infty} \lim_{N\to\infty} \, S_H^N(t) \, ,
$$
with 
\begin{equation} 
S_H^N(t) = \sum_{k=0}^H \frac{t^k}{k!} \int_{{\mathbb T}_N}
F_{0,N} \, L^k \varphi^{(m)} \, {\rm d} v_1\dots{\rm d} v_N\ ,
\label{eq_SHN(t)}
\end{equation}
Indeed, according to Lemma \ref{lem_absolute_convergence}, the series (\ref{eq_SHN(t)}) converges absolutely, uniformly with respect to $N$ and we can interchange the $H \to \infty$ and $N \to \infty$ limits. Now, for $k+m \leq N$, according to the inductive definition of $\varphi^{(m+k)}$ and to (\ref{eq_int_FN_Lphik}), we can write:
$$ \int_{{\mathbb T}_N}
F_{0,N} \, L^k \varphi^{(m)} \, {\rm d} v_1\dots{\rm d} v_N = \int_{{\mathbb T}_N}
F_{0,N} \,  \left( \varphi^{(m+k)} + \sum_{j=1}^k L^{k-j} \tilde \varphi^{(m+j)} \right) \, {\rm d} v_1\dots{\rm d} v_N . $$
Taking an index $H$ such that $H+m \leq N$, we have: 
\begin{eqnarray}
&& \hspace{-1cm}
S_H^N(t) = \sum_{k=0}^H \frac{t^k}{k!} \int_{{\mathbb T}_N}
F_{0,N} \,  \varphi^{(m+k)} \, {\rm d} v_1\dots{\rm d} v_N \, + \,  
\int_{{\mathbb T}_N}
F_{0,N} \, \,  R_H^N \varphi^{(m)} \, {\rm d} v_1\dots{\rm d} v_N , 
\label{eq_int_FN_etLphi}
\end{eqnarray}
with 
\begin{eqnarray*} 
R_H^N \varphi^{(m)} &=& \sum_{k=0}^H \frac{t^k}{k!}\, \sum_{j=1}^k L^{k-j} \tilde \varphi^{(m+j)} . 
\end{eqnarray*}
Now, taking successively $N \to \infty$ then $H \to \infty$, we show that the first term of (\ref{eq_int_FN_etLphi}) tends to the right-hand side of (\ref{eq_lim_int_FN_etLphi}) and the second one tends to zero. 

We start with the first term. By the chaos assumption, we have, as $N \to \infty$: 
\begin{eqnarray}
&& \hspace{-1cm} 
\lim_{N\to\infty} \sum_{k=0}^H \frac{t^k}{k!} \int_{{\mathbb T}_N}
F_{0,N} \, \varphi^{(m+k)} \, {\rm d} v_1\dots{\rm d} v_N = \nonumber \\
&& \hspace{3cm} 
=\sum_{k=0}^H \frac{t^k}{k!} \int_{{\mathbb T}_{m+k}}
\left(\prod_{j=1}^{k+m} f(v_j) \right) \, \varphi^{(m+k)} \, {\rm d} v_1\dots{\rm d} v_N 
. 
\label{eq_chaos_H}
\end{eqnarray}
We show that the series at the right-hand side of (\ref{eq_chaos_H}) is absolutely convergent, uniformly with respect to $t$ in any interval $[0,T]$ with $T<1/4$. The proof follows the same lines as that of Lemma \ref{lem_absolute_convergence} and we only sketch it. Using (\ref{induc}) and the Markov property of $Q_{(i,k+1)}$, we have 
\begin{eqnarray}
\| \varphi^{(m+k)} \|_\infty &\leq& 4 (m+k-1)  \, \| \varphi^{(m+k-1)} \|_\infty \nonumber \\
&\leq& 4^k \, k! \,  \left( \begin{array}{c} m+k-1 \\ k \end{array} \right) \,   \| \varphi^{(m)} \|_\infty \label{eq_phim+k} \\
&\leq& C_m 4^k \, k! \,  (k+1)^{m-1} \,   \| \varphi^{(m)} \|_\infty . \nonumber
\end{eqnarray}
Therefore, since $\prod_{j=1}^{k+m} f(v_j)$ is a probability:
\begin{eqnarray*}
&& \hspace{-1cm}
\left| \frac{t^k}{k!} \int_{{\mathbb T}_{m+k}}
\left(\prod_{j=1}^{k+m} f(v_j) \right) \, \varphi^{(m+k)} \, {\rm d} v_1\dots{\rm d} v_N  \right| \leq C_m (4t)^k \,  (k+1)^{m-1} \,   \| \varphi^{(m)} \|_\infty .
\end{eqnarray*}
This is the general term of a convergent series which converges uniformly as stated above. We deduce that 
\begin{eqnarray*}
&& \hspace{-1cm} 
\lim_{H\to\infty} \, \lim_{N\to\infty} \, \sum_{k=0}^H \frac{t^k}{k!} \int_{{\mathbb T}_N}
F_{0,N} \, \varphi^{(m+k)} \, {\rm d} v_1\dots{\rm d} v_N = \\
&& \hspace{4cm} 
=  \sum_{k=0}^\infty \frac{t^k}{k!} \int_{{\mathbb T}_{m+k}}
\left(\prod_{j=1}^{k+m} f(v_j) \right) \, \varphi^{(m+k)} \, {\rm d} v_1\dots{\rm d} v_N , 
\end{eqnarray*}
which is the right-hand side of (\ref{eq_lim_int_FN_etLphi}).

We now consider the second term of (\ref{eq_int_FN_etLphi}). First, using (\ref{eq_||Lkphi||_pvar_0}) with the pair $(k,p)$ replaced by $(k-j, m+j)$, we get 
\begin{eqnarray}
&& \hspace{-1cm} 
\| L^{k-j} \tilde \varphi^{(m+j)} \|_\infty \leq 4^{k-j} \, \, \frac{(m+k-1)!}{(m+j-1)!} \,  \,  \| \tilde \varphi^{(m+j)}\|_\infty \,
. 
\label{eq_Lk-jphim+j_1}
\end{eqnarray}
Then, (\ref{eq_estim_SN}) with $k$ replaced by $m+j-1$ yields 
\begin{eqnarray} 
&& \hspace{-1cm} 
\| \, \tilde \varphi^{(m+j)} \, \|_\infty \leq 6 \frac{(m+j-1)^2}{N-1}  \, \| \varphi^{(m+j-1)} \|_{\infty} \, . 
\label{eq_Lk-jphim+j_2}
\end{eqnarray}
Finally, using (\ref{eq_phim+k}) with $k=j-1$, we obtain
\begin{eqnarray}
\| \varphi^{(m+j-1)} \|_\infty \leq  4^{j-1} \, \frac{(m+j-2)!}{(m-1)!} \,   \| \varphi^{(m)} \|_\infty  . 
\label{eq_Lk-jphim+j_3}
\end{eqnarray}
Collecting (\ref{eq_Lk-jphim+j_1}), (\ref{eq_Lk-jphim+j_2}) and (\ref{eq_Lk-jphim+j_3}) leads to
\begin{eqnarray*}
&& \hspace{-1cm} 
\| L^{k-j} \tilde \varphi^{(m+j)} \|_\infty \leq \frac{3}{2} \, 4^k \, \, \frac{(m+k-1)!}{(m-1)!} \,  \, \frac{m+j-1}{N-1} \,  \| \varphi^{(m)}\|_\infty \,
, 
\end{eqnarray*}
and: 
\begin{eqnarray*}
\| \sum_{j=1}^k L^{k-j} \tilde \varphi^{(m+j)} \|_\infty & \leq & \frac{3}{2} \, 4^k \, \, \frac{(m+k-1)!}{(m-1)!} \,  \, \frac{(m+k)^2}{2(N-1)} \,  \| \varphi^{(m)}\|_\infty \, \\
& \leq & C_m \, 4^k \, k! \,  \, \frac{(k+1)^{m+1}}{N} \,  \| \varphi^{(m)}\|_\infty
, 
\end{eqnarray*}
which finally gives:
\begin{eqnarray*}
\frac{t^k}{k!} \| \sum_{j=1}^k L^{k-j} \tilde \varphi^{(m+j)} \|_\infty 
& \leq & C_m \, (4t)^k \,  \, \frac{(k+1)^{m+1}}{N} \,  \| \varphi^{(m)}\|_\infty
,
\end{eqnarray*}
Therefore, we have 
\begin{eqnarray*}
&&\hspace{-1cm}
\int_{{\mathbb T}_N}
F_{0,N} \, \,  R_H^N \varphi^{(m)} \, {\rm d} v_1\dots{\rm d} v_N \leq   \frac{1}{N} \, C_m \, \| \varphi^{(m)}\|_\infty \left( \sum_{k=0}^\infty (4t)^k \,  \, (k+1)^{m+1} \right)
,
\end{eqnarray*}
and deduce that 
\begin{eqnarray*}
\lim_{H\to\infty} \lim_{N\to\infty} \int_{{\mathbb T}_N}
F_{0,N} \, \,  R_H^N \varphi^{(m)} \, {\rm d} v_1\dots{\rm d} v_N = 0,
\end{eqnarray*}
which ends the proof. \endproof

\medskip
\noindent
We now have all the elements to prove Theorem \ref{thm_chaos_prop}. 

\medskip
\noindent
{\bf Proof of Theorem \ref{thm_chaos_prop}:} We need to show the existence of a function of a single velocity variable $f(v, t)$ such that 
\begin{eqnarray}
&& \hspace{-1cm} 
\lim_{N\to\infty} \int_{{\mathbb T}_N} F_{0,N} \, \, e^{tL} \varphi^{(m)} \, {\rm d} v_1\dots{\rm d} v_N 
= \int_{{\mathbb T}_m}
\left(\prod_{j=1}^{m} f(v_j,t) \right) \, \varphi^{(m)} \, {\rm d} v_1\dots{\rm d} v_N , 
\label{eq_etLFN_CP}
\end{eqnarray}
for all functions $\varphi^{(m)}$ of $m$ velocity variables $(v_1, \ldots, v_m)$ and all $m \leq 1$. 

First applying this definition with $m=1$ defines $f(v,t)$ (if it exists), by duality: 
\begin{eqnarray}
&& \hspace{-1cm} 
\int_{{\mathbb T}_1}
f(v,t) \, \varphi(v) \, {\rm d} v = \lim_{N\to\infty} 
\int_{{\mathbb T}_N} F_{0,N} \, \, e^{tL} \varphi \, {\rm d} v_1\dots{\rm d} v_N. 
\label{eq_def_f}
\end{eqnarray}
Now, applying (\ref{eq_etLFN_CP}) with $m=2$,
$\varphi^{(2)}(v_1,v_2) = \eta^{(1)}(v_1)\xi^{(1)}(v_2)$,
and using (\ref{eq_def_f})
leads to
\begin{multline}
\lim_{N\to\infty}\int_{{\mathbb T}_N} F_{0,N}\, e^{tL}\varphi^{(2)} \, {\rm d} v_1\dots{\rm d} v_N  =\\
\left(\lim_{N\to\infty}\int_{{\mathbb T}_N} F_{0,N} \, e^{tL}\eta^{(1)} \, {\rm d} v_1\dots{\rm d} v_N\right)
\left(\lim_{N\to\infty}\int_{{\mathbb T}_N} F_{0,N} \, e^{tL}\xi^{(1)} \, {\rm d} v_1\dots{\rm d} v_N\right) . 
\label{eq_CP_m=2}
\end{multline}
Reciprocally, (\ref{eq_CP_m=2}) implies (\ref{eq_etLFN_CP}) for $m=2$ and general  functions $\varphi^{(2)}(v_1,v_2)$ by the density of linear combinations of tensor products. Finally, the same argument applied to arbitrary $m$ would require to check the property for $\varphi^{(m)}$ equal to the tensor product of $m$ one-dimensional functions. This proof follows closely the proof for the case $m=2$ and will be omitted. So, we now focus to the case $m=2$. 

We first consider the case of small $t < 1/4$. Then, using (\ref{eq_lim_int_FN_etLphi}), eq. (\ref{eq_CP_m=2}) is equivalently written:
\begin{eqnarray*}
&& \hspace{-1cm}
\sum_{K=0}^\infty \frac{t^K}{K!} \int_{{\mathbb T}_{K+2}} \left(\prod_{j=1}^{K+2}f(v_j) \right) \varphi^{(K+2)} \, {\rm d} v_1\dots{\rm d} v_{K+2} = \\
&& \hspace{1cm}
\left( \sum_{k=0}^\infty \frac{t^k}{k!} \int_{{\mathbb T}_{k+1}} \left(\prod_{j=1}^{k+1}f(v_j) \right) \eta^{(k+1)} \, {\rm d} v_1\dots{\rm d} v_{k+1} \right) \times \\
&& \hspace{3cm}
\times
\left( \sum_{\ell=0}^\infty \frac{t^\ell}{\ell!} \int_{{\mathbb T}_{\ell+1}} \left(\prod_{j=1}^{\ell+1}f(v_j) \right) \eta^{(\ell+1)} \, {\rm d} v_1\dots{\rm d} v_{\ell+1} \right) \, .
\end{eqnarray*}
We start with the right-hand side of this formula, denoted by ${\mathcal R}$. By distributing the various terms in the product, we can write
\begin{eqnarray*}
&& \hspace{-1cm}
 {\mathcal R} = \sum_{K=0}^\infty \frac{t^K}{K!}  \int_{{\mathbb T}_{K+2}} \left(\prod_{j=1}^{K+2}f(v_j) \right) \left( \sum_{m=0}^K \left( \begin{array}{c} K \\ m \end{array} \right) \right.\\
&& \hspace{1cm}
\left. \phantom{\sum_{m=0}^K}  \eta^{(m+1)} (v_1, v_3, \ldots, v_{m+2}) \xi^{(K-m)+1} (v_2, v_{m+3}, \ldots, v_{K+2}) \right) \,  {\rm d} v_1 \dots {\rm d} v_{K+2} .
\end{eqnarray*}
Therefore, the result is proved if we show that for any symmetric function of $F(v_1, v_2,$ $\ldots, v_{K+2})$, we have 
\begin{eqnarray}
&& \hspace{-1cm}
\int_{{\mathbb T}_{K+2}} \left( \sum_{m=0}^K \left( \begin{array}{c} K \\ m \end{array} \right) \eta^{(m+1)} (v_1, v_3, \ldots, v_{m+2}) \xi^{(K-m)+1} (v_2, v_{m+3}, \ldots, v_{K+2}) \right) \nonumber \\
&& \hspace{7cm}
 \, F(v_1, v_2, \ldots, v_{K+2}) \,  {\rm d} v_1 \dots {\rm d} v_{K+2} = \nonumber \\
&& \hspace{0cm} 
= \int_{{\mathbb T}_{K+2}} \varphi^{(K+2)}(v_1, v_2, \ldots, v_{K+2}) \, F(v_1, v_2, \ldots, v_{K+2}) \,  {\rm d} v_1 \dots {\rm d} v_{K+2} 
 \, .
\label{eq_tensorization}
\end{eqnarray}
Indeed, at the initial step $K=1$, by direct computation we get:
\begin{eqnarray*}
&& \hspace{-1cm}
 \int_{{\mathbb T}_{3}} \varphi^{(3)}(v_1,v_2,v_3) \, F(v_1,v_2,v_3) \, {\rm d} v_1 \dots {\rm d} v_{3}=   \\
&& \hspace{1cm}
= \int_{{\mathbb T}_{3}} ( \eta^{(2)}(v_1,v_3)\xi^{(1)}(v_2) +   \eta^{(1)}(v_1)\xi^{(1)}(v_2,v_3)) \, F(v_1,v_2,v_3) \, {\rm d} v_1 \dots {\rm d} v_{3}  \ .
\end{eqnarray*}
Then, (\ref{eq_tensorization}) is easily proved by induction, using elementary properties of the binomial coefficients. This shows that $e^{t L^*} F_{0,N}$ is $f(\cdot, t)$-chaotic for small $t<1/4$. 

It remains to show that $f(v,t)$ is a solution of (\ref{eq_Boltzmann_PIM}). Again, assuming small $t<1/4$ and Using (\ref{eq_lim_int_FN_etLphi}) with $m=1$ and (\ref{eq_def_f}), we get 
\begin{eqnarray*}
&& \hspace{-1cm} 
\int_{{\mathbb T}_1}
f(v,t) \, \varphi^{(1)}(v) \, {\rm d} v = \sum_{k=0}^\infty \frac{t^k}{k!} \int_{{\mathbb T}_{k+1}} \left(\prod_{j=1}^{k+1}f(v_j) \right) \varphi^{(k+1)}(v_1,\dots,v_{k+1}) \, {\rm d} v_1\dots{\rm d} v_{k+1}\ .
\end{eqnarray*}
The convergence of the series at the right-hand side is uniform for $t<1/4$. So, we can differentiate this formula with respect to $t$ and obtain:
\begin{eqnarray}
&& \hspace{-1cm} 
\int_{{\mathbb T}_1}
\frac{\partial f}{\partial t} (v,t) \, \varphi^{(1)}(v) \, {\rm d} v = \nonumber \\
&& \hspace{1cm} 
= \sum_{k=0}^\infty \frac{t^k}{k!} \int_{{\mathbb T}_{k+2}} \left(\prod_{j=1}^{k+2}f(v_j) \right) \varphi^{(k+2)}(v_1,\dots,v_{k+2}) \, {\rm d} v_1\dots{\rm d} v_{k+2}\ .
\label{eq_partial_t_f}
\end{eqnarray}
Applying (\ref{eq_lim_int_FN_etLphi}) and (\ref{eq_etLFN_CP}) with $m=2$, we can re-write the right-hand side of (\ref{eq_partial_t_f}) and get
\begin{eqnarray*}
\int_{{\mathbb T}_1}
\frac{\partial f}{\partial t} (v,t) \, \varphi^{(1)}(v) \, {\rm d} v &=& 
\int_{{\mathbb T}_2}
\left(\prod_{j=1}^{2} f(v_j,t) \right) \, \varphi^{(2)} \, {\rm d} v_1 \, {\rm d} v_2 \\ 
&=& 
2 \int_{{\mathbb T}_2}
f(v_1,t) f(v_2,t) \, (Q_{(1,2)} - I) \varphi^{(1)} \, {\rm d} v_1 \, {\rm d} v_2 \\ 
&=& 
2 \int_{{\mathbb T}_2}
(Q_{(1,2)}^* - I)  (f \otimes f) \,  \varphi^{(1)} \, {\rm d} v_1 \, {\rm d} v_2 \ ,
\end{eqnarray*}
which is the weak form of (\ref{eq_Boltzmann_PIM}). This shows that $f(v,t)$ is a weak solution of (\ref{eq_Boltzmann_PIM}) for small $t<1/4$. 

In the last step of the proof, we need to remove the restriction on $t<1/4$. However, since this bound is independent of the initial data, we can partition any interval $[0,T]$ by intervals $[t_k, t_{k+1}]$ of length $t_{k+1} - t_k < 1/4$ and apply the result on each of these intervals with initial data $F_N(t_k)$ and $f(v,t_l)$. This shows that $e^{t L^*} F_{0,N}$ is $f(\cdot, t)$-chaotic on any finite-size interval $[0,T]$ where $f(v,t)$ is the solution of (\ref{eq_Boltzmann_PIM}). This ends the proof. \endproof

\begin{remark}
In Theorem \ref{thm_chaos_prop} the assumption that the pair selection probabilities $p_{i,j} (\vec v)$ are uniform, i.e. 
\begin{equation} 
p_{i,j} = \frac{2}{N(N-1)} , 
\label{eq_pair_uniform}
\end{equation}
is crucial. Indeed, we need at least that 

\medskip
\noindent (i) there exists a uniform constant  $C$ (independent of $N$) such that 
$$p_{i,j} \leq \frac{C}{N(N-1)}. $$
This is required to get the fundamental property (\ref{eq_fundamental_property}), with the constant $4k$ replaced by $4Ck$. 

\medskip
\noindent (ii) $p_{i,j} = p_{i,j}(v_i,v_j)$. 
In this way, the key fact in the proof of Lemma \ref{lem_absolute_convergence} that $L \varphi$, for $\varphi$ depending on only one velocity, is an average of terms only depending on two velocities and so on for $L^k \varphi$ remains true. 

\medskip 
\noindent (iii) $p_{i,j}(v,w) = p(v,w)$ is independent of $(i,j)$ to preserve the permutation symmetry of the problem. 

\medskip
\noindent
We see that these three properties together imply that (\ref{eq_pair_uniform}). Indeed, maintaining that 
\begin{equation} 
\sum_{i<j} p_{i,j} = 1, 
\label{eq_pair_normalization}
\end{equation}
for all velocity configurations necessitates that $p$ is a constant, and the normalization constraint (\ref{eq_pair_normalization}) leads to (\ref{eq_pair_uniform}). 
\label{rem_pair}
\end{remark}

\begin{remark}
Theorem \ref{thm_chaos_prop} can be extended to interaction processes involving multiple interactions as soon as the number of particles involved in an elementary interaction is finite and bounded independently of $N$. For instance, it will hold with a ternary interaction process, in which any interaction involves triples of particles. More generally, it will hold with a $p$-fold interaction process where any interaction involves exactly $p$ particles. As long as the interactions involve a finite number $p$ of interactions, with $p$ constant or bounded by a constant $P$ independent of $N$, the combinatorial arguments which have been developed above can be extended. Again, the master equation must combine the elementary interaction operators by means of uniform selection probabilities for the same arguments as those developed in Remark \ref{rem_pair}. Boltzmann operators with multiple interactions have been previously considered in \cite{BCG}. 
\label{rem_multiple} 
\end{remark}

\begin{remark}
As a by-product of Theorem \ref{thm_chaos_prop}, we have the existence of a weak solution of the nonlinear problem (\ref{eq_Boltzmann_PIM}). Additionally, the proof is constructive. However, the Theorem does not prove uniqueness, as another solution could exist not being related to a solution of the corresponding master equation. 
\label{rem_existence}
\end{remark}

\subsection{Application to the BDG and CL dynamics}
\label{subsec_CP_appli_BDG_CL}

\subsubsection{The BDG dynamics}
\label{subsecsec_CP_appli_BDG}

In the BDG case, thanks to (\ref{eq_def_Qij*_BDG}), we have
$$Q^*_{(1,2)} f^{\otimes 2} (v,w) =  \int_{{\mathbb T}_2} f(y_1) \, f(y_2) \, g(v\overline y_{i,2}^*) \, g(w\overline y_{1,2}^*) \, {\rm d} y_1 \, {\rm d} y_2\ ,$$
and thus
$$\int_{{\mathbb S}^1} Q^*_{(1,2)}f^{\otimes 2}(v,w) \, {\rm d} w = 
 \int_{{\mathbb T}_2} f(y_1) \, f(y_2) \, g(v\overline y_{i,2}^*) \, {\rm d} y_1 \, {\rm d} y_2\ .$$
Thus defining
$$Q_+(f,f)(v) =   \int_{{\mathbb T}_2} f(y_1) \, f(y_2) \, g(v\overline y_{i,2}^*) \, {\rm d} y_1 \, {\rm d} y_2\ ,$$
the one particle marginal $f(\cdot,t)$ such that  $\{e^{tL^*}F_{0,N}^{(N)}\}$ is $f(\cdot,t)$-chaotic
satisfies the kinetic-type equation
\begin{equation}\label{kke2}
\frac{\partial f}{\partial t}(v,t) = \int_{{\mathbb T}_2} f(y_1) \, f(y_2) \, g(v\overline y_{i,2}^*) \, {\rm d} y_1 \, {\rm d} y_2 - f(v,t)\ ,
\end{equation}
with $f(v,0) = f(v)$.

\subsubsection{The CL dynamics}
\label{subsecsec_CP_appli_CL}

In the CL case, thanks to (\ref{eq_def_Qij*_CL}) we have:
$$Q^*_{(1,2)}f^{\otimes 2}(v,w) =  \frac{1}{2}(f(v)+f(w)) \, g(v^*w)\ ,$$
and thus
$$\int_{\mathbb S}^1 Q^*_{(1,2)}f^{\otimes 2}(v,w) \, {\rm d} w = 
 \frac{1}{2} \left( f(v) + \int_{{\mathbb S}^1} f(w) \, g(v^*w) \, {\rm d} w \right)  = \frac{1}{2}[f(v) + f\star g(v)]\ ,$$
where $\star$ denote the convolution. 
Thus defining
$$Q_+(f,f)(v) =   \frac{1}{2}[f(v) + f\star g(v)]\ ,$$
the one particle marginal $f(\cdot,t)$ such that  $\{e^{tL^*}F_{0,N}\}$ is $f(\cdot,t)$-chaotic
satisfies the kinetic-type equation
\begin{equation}\label{kke3}
\frac{\partial f}{\partial t}(v,t)(v,t) = \frac{1}{2} [g\star f(v,t) - f(v,t)],
\end{equation}
with $f(v,0) = f(v)$. 

In this treatment, we have assumed that $g$ is independent of $N$. But if we let the variance of $g_N$ go to zero with $N$, we find
$$\lim_{N\to\infty}g_N\star f(v) = f(v)\ ,$$
and then we have
$$
\frac{\partial f}{\partial t}(v,t) = 0\ .
$$
That is, chaos is propagated, but nothing at all happens on the kinetic time scale. On a much longer time scale, correlations develop
and a new approach is needed to describe the bulk limit. This is what is shown in the next section by investigating the invariant densities, i.e. the equilibria of the master equation.

\setcounter{equation}{0}
\section{The Invariant densities $F_{\infty,N}$ for the CL dynamics}
\label{sec_invariant}

\subsection{Preliminaries}
\label{subsec_invariant_prliminaries}

Both the BDG and CL processes are clearly ergodic as long as $g$ is continuous, say, and so for  each  there will be a unique invariant density $F_\infty$, i.e., a unique density $F_\infty$ with $Q^*F_\infty = F_\infty$.  Since the process is symmetric under permutations of the variables, it is clear that
$F_\infty$ will be symmetric. It is not easy to write $F_\infty$ down in closed form. However, in the case of the CL dynamics, a very special property is true, namely the hierarchy of equations for the marginals (or BBGKY hierarchy) is closed at any order. This special feature will provide more information on the marginals of $F_\infty$. In particular, we show that under some specific scaling of the noise with respect to $N$, the invariant measure is not chaotic. This is not in contradiction to Theorem \ref{thm_chaos_prop}, since it is valid for fixed noise and on finite time intervals. 

The invariant density for the CL dynamics is the function $F_\infty$ which cancels the right-hand side of (\ref{eq_def_L*_BDG}) with $Q^*_{(i,j)}$ given by (\ref{eq_def_Qij*_CL}). Therefore, it satisfies
\begin{eqnarray}
&&\hspace{-1cm}
F_\infty(\vec v) = 
\frac{1}{N(N-1)} \sum_{i<j} \left[ 
[F_\infty]_{\widehat{i}}(v_1,\dots,\widehat{v_i},\dots,v_N)  + \right. \nonumber \\
&&\hspace{5cm} \left. +
 [F_\infty]_{\widehat{j}}(v_1,\dots,\widehat{v_j},\dots,v_N) \right] \, g(v_i^*v_j) \ .\label{update2}
\end{eqnarray}
While it is not easy to write $F_\infty$ down in closed form, we can at least say what $F_\infty$ is {\em not}: In general  $F(\vec v) =1$  {\em does not} solve $Q^*F = F$, i.e., $F_\infty$ is not the uniform density. Indeed, if we replace $F_\infty$ by $1$ on the right hand side of (\ref{update2}), we find
 $$\sum_{i<j}\frac{2}{N(N-1)}g(v_j^*v_i)\ .$$
This will equal $1$ for all $\vec v$ if and only if $g(z) =1$ for all $z$. 
However, it is easy to see that for fixed smooth $g$, and large N,
$$\sum_{i<j}\frac{2}{N(N-1)}g(v_j^*v_i) \approx 1\ ,$$
with high probability if $\vec v$ is selected at random, uniformly on ${\mathbb T}_N$, and so one can expect that for fixed $g$, the invariant density $F_\infty$ becomes more and more uniform as $N$ increases. 

However, a non-uniform invariant density can be found in the large $N$ limit if one includes some $N$ dependence in the noise density $g$ in such a way that it more and more closely approximates a $\delta$ function. This is perhaps justified in the context of biological modeling:  if the population is small, one fish may be less interested in carefully mimicking his neighbor than when the population is large.  One can imagine that the larger the group, the more important it is to follow behavioral rules closely. However, biological data are needed to support this claim.

While it does not seem easy to write $F$ down in closed form, it is possible to obtain analytical expressions of its marginals. This is the aim of the next section.

\subsection{Marginals}
\label{subsec_invariant_marginals}

For any symmetric density $F$ and any $m = 1,2,N-1$, define the $m$-variable marginal density
$F^{(m)}$ on ${\mathbb T}_m$ by
$$F^{(m)}(v_1,\dots,v_m) = \int_{{\mathbb T}_{N-m}}F(v_1,\dots,v_m,v_{m+1},\dots,v_N) \, {\rm d}v_{m+1}\dots {\rm d}v_N .$$
For typical binary collision process master equations, the evolution of the $m$-variable marginal density
depends on the $(m+1)$-variable marginal density, and one gets a hierarchy of evolution equations, the so-called BBGKY hierarchy. To break the hierarchy at some manageable level (e.g. $m=1$ in the usual kinetic theory case) in the large $N$ limit, one generally needs special assumptions on the initial data, such as the $f$-chaotic property, for a suitable $f$. One typically runs into the same hierarchy problem when trying to compute the marginals of the  invariant density, except that in some cases, one has $Q = Q^*$, and then the invariant density is simply uniform.

For the CL dynamics, the hierarchy breaks itself. Before stating the result, we recall that the Fourier series of a given function $\phi(w)$, $w \in {\mathbb S}^1$ is defined by
$$ \hat \phi(k) = \int_{{\mathbb S}^1} w^{-k} \, \phi(w) \, dw, \quad k \in {\mathbb Z}. $$
We have:  

\begin{proposition}
Let $F_\infty^{(1)}$ and $F_\infty^{(2)}$ be the one and two-variable marginal invariant densities of the CL process. We have:

\medskip
\noindent (i) $F_\infty^{(1)}$ is uniform, i.e. 
$$F_\infty^{(1)}(v_1)  =1, \quad \forall v_1 \in {\mathbb S}^1. $$

\medskip
\noindent (ii) $F_\infty^{(2)}$ is given by: 
$$F_\infty^{(2)}(v_1,v_2) = {\mathcal F}(v_1^*v_2), \quad \forall (v_1,v_2) \in {\mathbb T}_2 \ ,$$
with ${\mathcal F}(w)$, $w \in {\mathbb S}^1$ given by its Fourier series
\begin{eqnarray}
\widehat {\mathcal F}(k) 
&=& \frac{1}{N-1} \widehat g(k)\left[ 1 -  \frac{N-2}{N-1}\widehat g(k)\right]^{-1} ,  \label{form1} 
\end{eqnarray}
or equivalently 
\begin{eqnarray}
{\mathcal F} &=& \frac{1}{N-2}\sum_{\ell =1}^\infty \left[ \left(\frac{N-2}{N-1}\right)^\ell g^{*\ell}\right]\ .\label{form2}
\end{eqnarray}
\label{prop_marginals_invariant_CL}
\end{proposition}

\begin{remark}
Since $\sum_{\ell =1}^\infty  \left(\frac{N-2}{N-1}\right)^\ell = N-2$, 
Eq. (\ref{form2}) defines ${\mathcal F}$ as an average of convolution powers of $g$. 
\label{rem_mathcalF}
\end{remark}

\medskip
\noindent
{\bf Proof:} We begin with $F_\infty^{(1)}$. First of all, it is easy to see that for all $i  >1$ and all $j$ , then
$$
\int_{{\mathbb T}_{N-1}}  [F_\infty]_{\widehat{i}}(v_1,\dots,\widehat{v_i},\dots,v_N) g(v_j^*v_i) 
{\rm d}v_2\cdots  {\rm d}v_N =  F_\infty^{(1)}(v_1)\ ,
$$
Therefore, by integrating (\ref{update2}) with respect to $(v_2, \ldots, v_N)$ and using (\ref{eq_g_symmetric}), we have: 
\begin{eqnarray}
F_\infty^{(1)}(v_1) &=& 
\frac{1}{N(N-1)}\sum_{j=2}^N \int_{{\mathbb T}_{N-1}}
[F_\infty]_{\widehat{1}}(\widehat{v_1},\dots)   \, g(v_1^*v_j) \, {\rm d}v_2 \cdots  {\rm d}v_N
\nonumber\\
&&  \hspace{7.5cm} +  \frac{N-1}{N} \, F_\infty^{(1)}(v_1) . \label{lun2}
\end{eqnarray}
Then, using the permutation symmetry of $[F_\infty]_{\widehat{1}}$, one finds for each $j\ge 2$,
$$\int_{{\mathbb T}_{N-1}} [F_\infty]_{\widehat{1}}(\widehat{v_1},\dots) \, g(v_1^*v_j) \, {\rm d}v_2\cdots  {\rm d}v_N = 
\int_{{\mathbb S}^1} F_\infty^{(1)}(v_j) \, g(v_1^*v_j) \, {\rm d}v_j =  F_\infty^{(1)}*g(v_1)\ ,$$
where the $*$ denotes convolution on ${\mathbb S}^1$. Substituting it into (\ref{lun2}), we find
$$F_\infty^{(1)}(v_1) = \frac{1}{N}F_\infty^{(1)}*g(v_1) + \frac{N-1}{N}F_\infty^{(1)}(v_1)\ ,$$
which reduces to
$$F_\infty^{(1)}(v_1) =  F_\infty^{(1)}*g(v_1)\ .$$
The only solution of this equation for any $L^1$ density $g$ is the uniform density 
$F_\infty^{(1)}(v_1)  =1$.

We now turn to $F_\infty^{(2)}(v_1,v_2)$. Again, it is easy to see that if $i,j>2$, 
$$
\int_{{\mathbb T}_{N-2}}  [F_\infty]_{\widehat{i}}(v_1,v_2,\dots,\widehat{v_i},\dots,v_N) \, g(v_j^*v_i)  \, 
{\rm d}v_3\cdots  {\rm d}v_N =  F_\infty^{(2)}(v_1,v_2)\ ,
$$
and
$$
\int_{{\mathbb T}_{N-2}}  [F_\infty]_{\widehat{j}}(v_1,v_2,\dots,\widehat{v_j},\dots,v_N)  \, g(v_1^*v_j)  \, 
{\rm d}v_3\cdots  {\rm d}v_N =  F_\infty^{(2)}(v_1,v_2)\ .
$$
Next, for $i=1$ and $j=2$, using the permutation symmetry of $ [F_\infty]_{\widehat{1}}$ and (\ref{eq_g_symmetric}):
$$
\int_{{\mathbb T}_{N-2}}  [F_\infty]_{\widehat{1}}(\widehat{v_1},v_2,\dots,v_N)  \, g(v_2^*v_1)  \, 
{\rm d}v_3\cdots  {\rm d}v_N =  F_\infty^{(1)}(v_2) \, g(v_2^*v_1)  = g(v_1^*v_2)  \ ,
$$
and
$$
\int_{{\mathbb T}_{N-2}}  [F_\infty]_{\widehat{2}}(v_1,\widehat{v_2},\dots,v_N)  \, g(v_1^*v_2)  \, 
{\rm d}v_3\cdots  {\rm d}v_N =  F_\infty^{(1)}(v_1)g(v_1^*v_2)  = g(v_1^*v_2) \ .
$$
Finally, for $i =1$ and $j>2$, we have
$$
\int_{{\mathbb T}_{N-2}}  [F_\infty]_{\widehat{1}}(\widehat{v_1},v_2,\dots,v_N)  \, g(v_j^*v_1)  \, 
{\rm d}v_3\cdots  {\rm d}v_N = \int_{{\mathbb S}^1} F_\infty^{(2)}(v_2,v_j) \, g(v_j^*v_1) \, {\rm d}v_j   \ ,
$$
and
$$
\int_{{\mathbb T}_{N-2}}  [F_\infty]_{\widehat{j}}(v_1,v_2,\dots, \widehat{v_j},\dots,v_N)  \, g(v_1^*v_j)  \, 
{\rm d}v_3\cdots  {\rm d}v_N =  F_\infty^{(2)}(v_1,v_2) \ ,
$$
and for $i =2$ and $j>2$,  
$$
\int_{{\mathbb T}_{N-2}}  [F_\infty]_{\widehat{2}}(v_1,\widehat{v_2},\dots,v_N)  \, g(v_j^*v_2)  \, 
{\rm d}v_3\cdots  {\rm d}v_N = \int_{{\mathbb S}^1} F_\infty^{(2)}(v_1,v_j) \, g(v_j^*v_2) \, {\rm d}v_j   \ ,
$$
and
$$
\int_{{\mathbb T}_{N-2}}  [F_\infty]_{\widehat{j}}(v_1,v_2,\dots, \widehat{v_j},\dots,v_N)  \, g(v_2^*v_j)  \, 
{\rm d}v_3\cdots  {\rm d}v_N =  F_\infty^{(2)}(v_1,v_2) \ .
$$

We can now compute the two variable marginals of both sides of (\ref{update2}), and we find
\begin{eqnarray*}
&&\hspace{-1cm}
F_\infty^{(2)}(v_1,v_2)  = 
\left[\frac{(N-2)(N-3)}{N(N-1)} + 2\frac{N-2}{N(N-1)}\right] F_\infty^{(2)}(v_1,v_2)  \\
&&\hspace{4cm}
+\frac{2}{N(N-1)}g(v_1^*v_2) 
+ \frac{N-2}{N(N-1)} H(v_1,v_2),
\end{eqnarray*}
where
$$H(v_1,v_2) =  \int_{{\mathbb S}^1} F_\infty^{(2)}(v_2,z) \, g(z^*v_1) \, {\rm d}z+
\int_{{\mathbb S}^1} F_\infty^{(2)}(v_j,v_1) \, g(z^*v_2) \, {\rm d}z\ .$$
This simplifies to 
$$F_\infty^{(2)}(v_1,v_2)  = \frac{1}{N-1}g(v_1^*v_2) + \frac{N-2}{2(N-1)}H(v_1,v_2) \ .$$
 
Fourier transforming both sides, we have
$$\widehat H(k_1,k_2) =  \widehat F_\infty^{(2)}(k_1,k_2)(\widehat g(k_1) +  \widehat g(k_2) )$$
and so
$$\widehat F_\infty^{(2)}(k_1,k_2)\left[ 1 - \frac{N-2}{N-1}\left(\frac{\widehat g(k_1) +  \widehat g(k_2)}{2}\right)\right] 
 = \frac{1}{N-1}\widehat g(k_1) \delta_{k_1,-k_2}\ ,$$
where $\delta_{i,j}$ is the usual Kronecker symbol. 
It follows that  $\widehat F_\infty^{(2)}(k_1,k_2)$ has the form $\widehat {\mathcal F}(k_1) \delta_{k_1,-k_2}$, i.e. that $F_\infty^{(2)}$ has the form $F_\infty^{(2)}(v_1,v_2) = {\mathcal F}(v_1^*v_2)$ with ${\mathcal F}$ defined by its Fourier transform according to (\ref{form1}) (owing to the evenness of $\hat g$, a consequence of (\ref{eq_g_symmetric})). This ends the proof of proposition \ref{prop_marginals_invariant_CL}. \endproof

\subsection{Noise scaling}
\label{subsec_invariant_noise}

The reason for considering a scaling of the noise intensity is the following. For fixed $g\in L^1({\mathbb S}^1)$, $\lim_{\ell\to\infty}g^{*\ell} =1$, and since for large $N$,
most of the weight in the average is on large values of $\ell$, ${\mathcal F}$ will be nearly uniform for large values of $N$, and the correlations are washed out. Therefore, we recover here that the invariant density is nearly uniform for large values of $N$, a fact which has already been noticed (see section \ref{subsec_invariant_prliminaries}). 
 
But if $g$ is taken to depend on $N$ itself, this need not be the case. As a typical example, we can consider $g(z^*w)$ to be the kernel of $e^{\Delta/N}$ on ${\mathbb S}^1$; i.e., the heat kernel on ${\mathbb S}^1$ at time $1/N$.
Then 
\begin{equation}
\widehat g_N(k) := e^{-k^2/N}\ .
\label{eq_heat}
\end{equation}
We now state the 

\begin{proposition}
Suppose the scaled noise intensity $g_N$ is such that 
\begin{equation}
\lim_{N\to\infty} (N-2)(\widehat g_N(k)-1) := \gamma(k)
\label{eq_def_gamma}
\end{equation}
exists and is non trivial (i.e. not equal to the Kronecker $\delta_{k,0}$). Then, the corresponding correlation ${\mathcal F}_N$ associated to $g_N$ through (\ref{form1}) or (\ref{form2}) satisfies:
$$ \lim_{N\to\infty} \widehat {\mathcal F}_N (k) := {\mathcal F}_\infty (k) = \frac{1}{1 - \gamma(k)}. $$
\label{prop_noise_scaling}
\end{proposition}

\medskip
\noindent
{\bf Proof:} Using $g_N$  in place of $g$ in (\ref{form1}), we find
\begin{eqnarray}\label{kj}
 \widehat {\mathcal F}_N(k) &= &\frac{1}{N-1}\widehat g_N(k)\left[ 1 -  \frac{N-2}{N-1}(1 + (\widehat g_N(k)-1))\right]^{-1}\nonumber\\
  &= &\frac{1}{N-1}\widehat g_N(k)\left[  \frac{1}{N-1} -  \frac{N-2}{N-1}(\widehat g_N(k)-1)\right]^{-1}\nonumber\\
 &= &\widehat g_N(k)\left[ 1 -  (N-2)(\widehat g_N(k)-1)\right]^{-1}\ .\nonumber\\
  \end{eqnarray}
The result follows from inserting (\ref{eq_def_gamma}). \endproof

\begin{example} The heat kernel (\ref{eq_heat}). In this case, we have 
$$ \gamma(k ) = - k^2, \qquad {\mathcal F}_\infty (k) = \frac{1}{1 + k^2}. $$
Hence, the correlation function is a Lorentzian in Fourier space. 
\end{example}

With this noise scaling, the two-variable marginal invariant density of the $N$-particle CL process $F_{\infty,N}^{(2)}(v_1,v_2) $ is such that 
$$ F_{\infty,N}^{(2)}(v_1,v_2) \to {\mathcal F}_\infty (v_1^*v_2), \quad \mbox{ as } \quad N \to \infty, $$
where ${\mathcal F}_\infty$ is not the uniform distribution. Therefore, non-trivial correlations remain in the large $N$ limit and in particular, the invariant density  $\{F_{\infty,N}\}$ is not chaotic. This result is in marked contrast to the case studied by Kac, in which
the invariant density is the uniform density on the sphere $S^{N-1}(\sqrt{N})$.  This family is well known to be $G$-chaotic where
$G(v)$ denotes the centered unit Gaussian on ${\mathord{\mathbb R}}$.  The lack of chaos in the invariant density might seem to be a strong obstacle to propagation of chaos. But as we have seen in Theorem \ref{thm_chaos_prop}, this is not the case. As already noticed, there is no contradiction between this two seemingly paradoxical results. Theorem \ref{thm_chaos_prop} is valid under fixed noise and on a finite time interval. By contrast, the lack of chaos property of the invariant density for the CL model is shown under $N$-dependent noise intensity and in the infinite time limit. But these two results show that the chaos property can be valid for a finite time interval at the kinetic time scale and be lost at larger times.

\setcounter{equation}{0}
\section{Conclusion}
\label{sec_conclu}

We have considered a class of pair-interaction stochastic processes in an $N$-particle system and their associated pair interaction driven master equations. We have proved a chaos propagation theorem for this class of master equations and have used this result to study the kinetic limits of two biological swarm models, the BDG and CL processes. By investigating the invariant density of the CL process, we have shown that the chaos property may be lost at large times. This work shows that the chaos property may be true even for processes that seemingly build-up correlations but may not be uniformly valid in time. Correlation build-up manifests itself at large time scales. In order to restore the validity of kinetic theory at these large scales, new theories must must be developed. This is a fascinating and widely open area of research.


\end{document}